\documentclass[19pt]{article}
\usepackage[utf8]{inputenc}
\usepackage[english,russian]{babel}
\usepackage{indentfirst}
\usepackage{misccorr}
\usepackage{color}
 \usepackage{amsmath}
\usepackage{amssymb}

\def\R{{\mathbb R}}
\def\aff{\mbox{aff}}

\begin{document}

\title{A `converse' to the Constraint Lemma}

\author{E. Kolpakov\thanks{Supported by the Russian Foundation
for Basic Research grant 19-01-00169}}

\date{}

\maketitle

We begin the article with the statement of the main result.
For the history and motivation see Remark 1; the current paper gives a partial answer to [Sk18, Problem 4.4.c].
Direct proofs of implications between classical theorems are interesting.
The main result is a direct proof of the implication $(LVKF_{k,3})\Rightarrow( LTT_{3k-1,3})$ below.
Consider the following statements:

\smallskip
($LVKF_{1,3}$) From any 11 points in $ \mathbb{R}^{3}$ one can choose 3 pairwise disjoint triples whose convex hulls have a common point.

($LVKF_{k,3}$) From any $6k + 5$ points in $ \mathbb{R}^{3k}$ one can choose 3 pairwise disjoint sets each containing  $2k + 1 $ points and whose convex hulls have a common point.

($LTT_{2,3}$)  Any 7 points in $\mathbb{R}^{2}$  can be decomposed into 3 subsets whose convex hulls have a common point.

($LTT_{d,3}$)  Any $2d+3$ points in $\mathbb{R}^d$  can be decomposed into 3 subsets whose convex hulls have a common point.

\smallskip
These statements are correct and known, but the purpose of this paper is the direct deduction of one statement from another.

\smallskip
{\bf  Remark 1.}
Consider the following piecewise linear versions (or, equivalently, topological versions) of the above statements.
Denote by $\Delta_{N}$ the $N$-dimensional simplex.

($VKF_{k, r}$)  (the $r$-fold van Kampen-Flores conjecture) {\it For any general position PL (piecewise linear) map $\Delta_{(kr+2)(r-1)}\to\R^{kr}$ there are $r$ pairwise disjoint $k(r-1)$-faces whose images have a common point.}

($TT_{d, r}$)  (the topological Tverberg conjecture) {\it For any PL map $\Delta_{(d+1)(r-1)}\to \R^d$ there are $r$ pairwise disjoint faces whose images have a common point.}

These conjectures are true for $r$ a prime power, and are false for $r$ not a prime power, see the surveys [BZ], [BBZ], [Sk18].

The implication $(TT_{kr+1, r})  \Rightarrow  (VKF_{k, r})$ is called the Constraint Lemma [Gr10, BFZ14, Fr15], for a history see [Sk18, Remark 1.11].
This is one of the steps (the simplest) in the disproof of the topological Tverberg conjecture, see the surveys [BZ], [BBZ], [Sk18].

Our result $(LVKF_{k,3})\Rightarrow( LTT_{3k-1,3})$ can be considered as a `converse' to the Constraint Lemma.
However, $( LTT_{d,3})$ is correct for every $d$, so our direct proof does not give any new results.

In addition to the above, there are the following implications for $r=2$, see the survey [Sk17]:

$TT_{d, 2}  \Rightarrow   TT_{d-1, 2}$;

$ TT_{d, 2}  \Rightarrow  VKF_{d-1, 2}$ for even $d$;

$VKF_{d, 2}  \Rightarrow   VKF_{d-1, 2}$;

$VKF_{d-1, 2}  \Rightarrow   VKF_{d, 2}$ for odd $d$;


For direct proofs of other similar implications see [Sk17] and references therein.
\smallskip



The most novel part of our proof a simple proof of lemma 1, which is a weaker version of [Tv66, Lemma 3].

Lemmas 2,3,4 below are simple and so could be known.

Denote the convex hull of a set $ P $ by $ \left< P \right> $.
A set of points in the $\R^n$ is called {\it in general position} if it contains no $k$ points in $(k-2)$-dimensional subspace for any $k\leq n+1$.

\smallskip
{\bf Lemma 1.} {\it Let

$\bullet$ $A$ be a set of $6k+1$ points of general position in $\R^{3k-1}\times 0\subset \R^{3k}$;

$\bullet$ $A_1$ be a vertex of the convex hull  $\left<A\right>$;


$\bullet$ $A_1\times0,M_1, M_2, M_3, M_4$ be points in the line  $A_1\times \R\subset \R^{3k-1}\times\R=\R^{3k}$ in order of increasing the last coordinate;


$\bullet$ $\Delta_1, \Delta_2, \Delta_3 \subset A\cup\{ M_1, M_2, M_3, M_4\}$ be pairwise disjoint subsets such that $M_2\notin \left<\Delta_1\right>\cap \left<\Delta_2\right>\cap \left<\Delta_3\right> \neq \varnothing$.

Then there exist pairwise disjoint subsets  $\Delta'_1, \Delta'_2, \Delta'_3 \subset A\cup\{M_1, M_2, M_3, M_4\}$ such that
\begin{enumerate}
  \item $\left<\Delta'_1\right>\cap \left<\Delta'_2\right>\cap \left<\Delta'_3\right> \neq \varnothing$; 
  \item  some two the points of $A_1, M_1, M_2, M_3, M_4$ do not belong to 
  $\Delta'_1\cup\Delta'_2\cup\Delta'_3$; 
  \item $(\Delta'_1\cup\Delta'_2\cup\Delta'_3)\cap A_1\times \R \subset(\Delta_1\cup\Delta_2\cup\Delta_3)\cap A_1\times \R$.
\end{enumerate}
}

\smallskip
{\it Proof}.
  The point $X$ is called {\it higher} than $Y$ if last coordinate of $X$ is greater than last coordinate of $Y$.

 Consider all triples of  pairwise disjoint subsets $\delta_1, \delta_2, \delta_3 \subset A\cup\{M_1, M_2, M_3, M_4\}$ such that  $(\delta_1\cup \delta_2\cup \delta_3)\cap A_1\times \R \subset(\Delta_1\cup\Delta_2\cup\Delta_3)\cap A_1\times \R$ . The union of the intersections  $ \left<\delta_1\right>\cap \left<\delta_2\right>\cap \left<\delta_3\right>$ is non-empty, closed and bounded. Denote by $Z'$ the point with the smallest last coordinate from this union.
Therefore $Z'\in \left<\Delta'_1\right>\cap \left<\Delta'_2\right>\cap \left<\Delta'_3\right>$ for some pairwise disjoint subsets $\Delta'_1, \Delta'_2, \Delta'_3 \subset A\cup\{ M_1, M_2, M_3, M_4\}$. From the sets $\Delta'_1, \Delta'_2, \Delta'_3$ there is at least one, which not contain the points $M_3, M_4$. This set is contained in $A\cup\{M_1,M_2\}$, and its convex hull besides the point $M_2$ is below than the point $M_2$. Hence the point $Z'$ is below than the point $M_2$. Since the point $A_1$ is a vertex of $\left<A\right>$, at least one of the sets $\left<\Delta'_1\right>, \left<\Delta'_2\right>, \left<\Delta'_3\right>$ intersects the line $A_1\times \R$ in no more than one point. Since $Z'\neq M_2$,  $Z'$ does not belong to line $A_1\times \R$.

  Each of polyhedrons $\left<\Delta'_1\right>, \left<\Delta'_2\right>, \left<\Delta'_3\right>$ can be replaced by a simplex from its triangulation containing the point $Z'$. Therefore  without loss of generality we consider that the polyhedrons $\left<\Delta'_1\right>, \left<\Delta'_2\right>, \left<\Delta'_3\right>$ are simplexes. Replace each simplexes $\left<\Delta'_1\right>, \left<\Delta'_2\right>, \left<\Delta'_3\right>$ by its face such that the point $Z'$ became is unique point of the intersection $\left<\Delta'_1\right>\cap \left<\Delta'_2\right>\cap \left<\Delta'_3\right>$, and does not belong the faces of the simplexes $\left<\Delta'_1\right>, \left<\Delta'_2\right>, \left<\Delta'_3\right>$.  Since  $\dim\left<\Delta'_1\right>+\dim\left<\Delta'_2\right>+\dim\left<\Delta'_3\right>\leq 6k$, we have $|\Delta'_1|+|\Delta'_2|+|\Delta'_3|\leq 6k+3$.

 If any point from $A-\{A_1\}$ belong to some set from $\Delta'_1, \Delta'_2, \Delta'_3$ then from the points $A_1, M_1, M_2, M_3, M_4$ at least two do not belong in any of the sets $\Delta'_1, \Delta'_2, \Delta'_3$ as required.

 If $Z'$ belong to plane $\R^{3k-1} \times 0$ then for any $i=1,2,3$ we have $Z'\in\R^{3k-1} \times 0 \cap\left<\Delta'_i\right>=\left<\alpha\cap\Delta'_i\right>$. Therefore the sets $\R^{3k-1} \times 0 \cap\Delta'_1,\R^{3k-1} \times 0\cap\Delta'_2,\R^{3k-1} \times 0\cap\Delta'_3$ are required.

 Next, we consider the remaining case where exists a point $Y$ from the set $A-\{A_1\}$ that does not belong to any of the sets  $\Delta'_1, \Delta'_2, \Delta'_3$, and the point $Z'$ does not lie on the plane $\R^{3k-1} \times 0$. Then $Z'$ have positive last coordinate. Consider a convex polyhedron $\mathfrak{M}:=\left<Y\cup\Delta'_1\right>\cap \left<Y\cup\Delta'_2\right>\cap \left<Y\cup\Delta'_3\right>$, and proove that $Z'$ is a his vertex. Suppose that $Z'$ is not a vertex of  polyhedron $\mathfrak{M}$, then there is a segment  $I\subset\mathfrak{M}$, whose middle point is $Z'$. For any $i=1,2,3$ we have: $I\subset\left<Y\cup\Delta'_i\right>$, $Z'$ does not belong to face of simplex $\left<\Delta'_i\right>$. By Lemma 2 for any $i=1,2,3$ there is a segment $J_i\subset I$ such that $J_i\subset\left<\Delta'_i\right>$ and $Z'$is a middle $J_i$. Then for $J=J_1\cap J_2\cap J_3$ we have $J\subset\left<\Delta'_1\right>\cap \left<\Delta'_2\right>\cap \left<\Delta'_3\right>$ but $Z'$ is a unique point of the intersection  $\left<\Delta'_1\right>\cap \left<\Delta'_2\right>\cap \left<\Delta'_3\right>$. This contradiction means that  $Z'$  is the vertex of the polyhedron  $\mathfrak{M}$.

 For any $i=1,2,3$ denote $\tau_i=\Delta'_i\cap(A-\{A_1\})$. Denote by $Z''$ orthogonal projection of point $Z'$ onto $\R^{3k-1} \times 0$. Then $Z''\in \left<A_1\cup\tau_1\right>\cap\left<A_1\cup\tau_2\right>\cap\left<A_1\cup\tau_3\right>$ and $Z''\neq A_1$. Consider the ray $L$ with vertex in $A_1$ which contains a point  $Z''$. This ray intersect a simplexes $\left<\tau_1\right>,\left<\tau_2\right>,\left<\tau_3\right>$. For any $i=1,2,3$ denote the point of intersection (possibly non-single) of the ray $L$ and simplex $\left<\tau_i\right>$ by $V_i$. The segments $[A_1;V_1], [A_1;V_2], [A_1;V_3]$ are ordered by inclusion. Without less of generality let $[A_1;V_3]\subset [A_1;V_2]\subset [A_1;V_1]$.
Then $V_3\in\left<A_1\cup\tau_1\right>\cap\left<A_1\cup\tau_2\right>\cap\left<\tau_3\right>$ and $\left<A_1\cup\tau_1\right>\cap\left<A_1\cup\tau_2\right>\cap\left<\tau_3\right> \neq\varnothing$.

  In the next two paragraphs prove that $|\Delta'_1|+|\Delta'_2|+|\Delta'_3|= 6k+3$.

 If there is a simplex among  $\left<\Delta'_1\right>, \left<\Delta'_2\right>, \left<\Delta'_3\right>$ which dimension is  $3k$ then it contains the point  $Z'$ as internal.
Denote by $(i_1, i_2, i_3)$ a permutation of numbers $(1, 2, 3)$ such that $\dim\left<\Delta'_{i_3}\right>=3k$ and $i_1< i_2$
Since the set $A$ is in general position in  $\R^{3k-1} \times 0$, we have $\dim\left<A_1\cup\tau_{i_1}\right>\cap \left<A_1\cup\tau_{i_2}\right>= 1+\dim\left<A_1\cup\tau_{i_1}\right>\cap \left<\tau_{i_2}\right>=|\tau_{i_1}|+|\tau_{i_2}|-(3k-1)$.
Since $[A_1;Z'']\subset\left<A_1\cup\tau_{i_1}\right>\cap \left<A_1\cup\tau_{i_2}\right>$, we have $\dim\left<A_1\cup\tau_{i_1}\right>\cap \left<A_1\cup\tau_{i_2}\right>\geq 1$. Hence $|\tau_{i_1}|+|\tau_{i_2}|\geq 3k$ therefore
$|\Delta'_{i_1}|+|\Delta'_{i_2}|+|\Delta'_{i_3}|\geq (|\tau_{i_1}|+1)+(|\tau_{i_2}|+1)+(3k+1)\geq 6k+3$ as required.

 If dimensions of simplexes $\left<\Delta'_1\right>, \left<\Delta'_2\right>, \left<\Delta'_3\right>$ not more than $3k-1$ then dimensions of simplexes $\left<\tau_1\right>, \left<\tau_2\right>, \left<\tau_3\right>$ not more than $3k-2$
and this simplexes are faces of simplexes $\left<A_1\cup\tau_1\right>, \left<A_1\cup\tau_2\right>, \left<A_1\cup\tau_3\right>$ respectively.
Since the points of the set $\{A_1\}\cup\tau_1\cup\tau_2$ is in general position, we have $\dim\left<A_1\cup\tau_1\right>\cap \left<A_1\cup\tau_2\right>= 1+\dim\left<\tau_1\right>\cap \left<A_1\cup\tau_2\right>=|\tau_1|+|\tau_2|-(3k-1)$.
Since the points of the set $\{A_1\}\cup\tau_1\cup\tau_2$  is in general position with the points of the set $\tau_3$, we have
$$\dim\left<A_1\cup\tau_1\right>\cap\left<A_1\cup\tau_2\right>\cap\left<\tau_3\right>=$$
$$=\dim\left<A_1\cup\tau_1\right>\cap\left<A_1\cup\tau_2\right>+\dim\left<\tau_3\right>-(3k-1) =$$
$$=|\tau_1|+|\tau_2|-(3k-1)+(|\tau_3|-1)-(3k-1)=|\tau_1|+|\tau_2|+|\tau_3|-(6k-1)$$
Since $\left<A_1\cup\tau_1\right>\cap\left<A_1\cup\tau_2\right>\cap\left<\tau_3\right> \neq\varnothing$, we have $|\tau_1|+|\tau_2|+|\tau_3|-(6k-1)=\dim\left<A_1\cup\tau_1\right>\cap\left<A_1\cup\tau_2\right>\cap\left<\tau_3\right>\geq 0$.
Since from the points $A_1, M_1,  M_2,  M_3,  M_4$ at least 4 belong to the union $\Delta'_1\cup\Delta'_2\cup\Delta'_3$, we have $|\Delta'_1|+|\Delta'_2|+|\Delta'_3|\geq 4+|\tau_1|+|\tau_2|+|\tau_3|= 6k+3$ as required.

 Since $|\Delta'_1|+|\Delta'_2|+|\Delta'_3|= 6k+3$, we have $|Y\cup\Delta'_1|+|\Delta'_2|+|\Delta'_3|= 6k+4$. Therefore
$\dim\left<Y\cup\Delta'_1\right>+\dim \left<\Delta'_2\right>+\dim \left<\Delta'_3\right>\geq 6k+1$.
Hence $\dim\left<Y\cup\Delta'_1\right>\cap \left<\Delta'_2\right>\cap \left<\Delta'_3\right>\geq 1$.
It is similarly proved that $\dim\left<\Delta'_1\right>\cap \left<Y\cup\Delta'_2\right>\cap \left<\Delta'_3\right>\geq 1$ and
$\dim\left<\Delta'_1\right>\cap \left<\Delta'_2\right>\cap \left<Y\cup\Delta'_3\right>\geq 1$.

 The pairwise intersections of the sets $\left<Y\cup\Delta'_{i_1}\right>\cap \left<\Delta'_{i_2}\right>\cap \left<\Delta'_{i_3}\right>$ for all permutations $\{i_1, i_2, i_3\}=\{1,2,3\}$ are an $\left<\Delta'_1\right>\cap \left<\Delta'_2\right>\cap \left<\Delta'_3\right>=Z'$. Since $\left<Y\cup\Delta'_{i_1}\right>\cap \left<Y\cup\Delta'_{i_2}\right>\cap \left<\Delta'_{i_3}\right>$ contains $\left<Y\cup\Delta'_{i_1}\right>\cap \left<\Delta'_{i_2}\right>\cap \left<\Delta'_{i_3}\right>$ and $\left<\Delta'_{i_1}\right>\cap \left<Y\cup\Delta'_{i_2}\right>\cap \left<\Delta'_{i_3}\right>$, we have $\dim\left<Y\cup\Delta'_{i_1}\right>\cap \left<Y\cup\Delta'_{i_2}\right>\cap \left<\Delta'_{i_3}\right>\geq 2$ for $\{i_1, i_2, i_3\}=\{1,2,3\}$. By Lemma 3  any face of $\left<Y\cup\Delta'_{i_1}\right>\cap \left<Y\cup\Delta'_{i_2}\right>\cap \left<\Delta'_{i_3}\right>$ which contains the point $Z'$, is contained  in some of the sets $\left<Y\cup\Delta'_{i_1}\right>\cap \left<\Delta'_{i_2}\right>\cap \left<\Delta'_{i_3}\right>$ or $\left<\Delta'_{i_1}\right>\cap \left<Y\cup\Delta'_{i_2}\right>\cap \left<\Delta'_{i_3}\right>$.

 Then any edge of the polyhedron $\mathfrak{M}$ with vertex in $Z'$ is contained in one of the sets $\left<Y\cup\Delta'_1\right>\cap \left<\Delta'_2\right>\cap \left<\Delta'_3\right>$, $\left<\Delta'_1\right>\cap \left<Y\cup\Delta'_2\right>\cap \left<\Delta'_3\right>$, $\left<\Delta'_1\right>\cap \left<\Delta'_2\right>\cap \left<Y\cup\Delta'_3\right>$. Since the point $Z'$ above than the point  $Y$,  $Z'$ is not the lowest point of a polyhedron $\mathfrak{M}$. Then $Z'$ is not the lowest vertex of a polyhedron $\mathfrak{M}$. Therefore there is edge $[Z';Z_1]$ such that the point $Z_1$ is an below $Z'$.  Without loss of generality  $[Z';Z_1]\subset\left<Y\cup\Delta'_1\right>\cap \left<\Delta'_2\right>\cap \left<\Delta'_3\right>$. Since the point $Z_1$ is an below $Z'$, we have contradiction with choice of the point $Z'$.
 Lemma is proved.

\smallskip

{\bf Lemma 2.} {\it Take a simplex $S$ in the space $\R^n$.
Take a point $Y$ in general position with vertexes of $S$.
Take a segment $I\subset\left<Y\cup S\right>$.
Denote by $Z$ the middle point of $I$.
Suppose that $Z\in S$ but $Z$ does not belong any face of $S$.
Then there is a non-degenerate segment $J\subset I$  such that $J\subset S$ and $Z$ is a midpoint of $J$}.

\smallskip
{\it Proof}.
 If $\dim S=n$, then $Z$ is interior point of simplex $S$. There is closed $\epsilon$-neighborhood   of the point $Z$ which belong $S$ for some $\epsilon>0$. Take as  $J$ the intersection of such a neighborhood with  $I$. Then $Z$ will be the middle of the segment  $J$ and $J\subset S$ as required.

 If $\dim S\leq n-1$, then $S$ is a face of the simplex $\left<Y\cup S\right>$. Since the point $Z$ is a middle of the segment $I$ and  $Z\in S$, we have that both ends of the segment $I$ belong to $S$. Take $J=I$. Lemma is proved.

\smallskip

{\bf Lemma 3.} {\it Take convex polyhedra $A_1, ... , A_n$ in $\R^d$ which have more than one common point.
Then for any point $P$ in any face of  $A_1\cap ... \cap A_n$ there are faces $A_j'$ of $A_j$, $j=1,2,\ldots,n$,  such that}
$$P\in (A'_1\cap A_2\cap ... \cap A_n)\cup(A_1\cap A_2'\cap A_3\cap ... \cap A_n)\cup\ldots\cup (A_1\cap A_2\cap ... \cap A_n').$$ 

 {\it Proof}. It suffices to prove the proposition for $n=2$.

 Suppose that there is a point $P\in (A_1\cap A_2)'$ such that $P \notin A'_1$ and $P\notin A'_2$.
 Denote by $d_1$ and $d_2$ the dimensions of the polyhedra $A_1$ and $A_2$ respectively. Take the disks $D_1$ and $D_2$  with dimensions $d_1$ and $d_2$ respectively such that $P\in D_1\subset A_1$ and $P\in D_2\subset A_2$. The intersection $D_1\cap D_2$ is a disk with dimension $d_1+d_2-d=\dim A_1\cap A_2$. Since $D_1\cap D_2\subset A_1\cap A_2$, the point $P$ does not belong to face of $A_1\cap A_2$. A contradiction.
 
\smallskip

Denote the affine hull of a set $ P $ by $\aff  P$.

{\bf Lemma 4.} {\it Let   $P$ a finite general position subset of $ \mathbb{R}^{3k-1}$  and $\epsilon>0$ a number. Then there is $\delta \in (0,\epsilon)$ with the following property. For any two subsets  $P_1,P_2$ of $P$  such that the intersection $\left<P_1\right>\cap \left<P_2\right>$ is one point
we have that $\epsilon$-neighbourhood  of this point $\left<P_1\right>\cap \left<P_2\right>$ contains the intersection of $\delta$-neighbourhoods of
$\left<P_1\right>$ and $\left<P_2\right>$.}

\smallskip
{\it Proof}.
Since $P$ is finite, there is only a finite number of pairs $(P_1,P_2)$ of subsets of $P$ such that $\left<P_1\right>\cap \left<P_2\right>$
is a point, but none of the subsets is a point.
Let $\alpha>0$ be the smallest of any angles  between $\left<P_1\right>$ and $\left<P_2\right>$.
(The {\it angle} between two intersecting affine subspaces of $\R^{3k-1}$,  none of which is a point, is the minimal
angle beetwen two rays originating from an intersection point of the subspaces,
the first ray contained in the first subspace and the second ray in the second subspaces.)
Let $\delta=\epsilon \sin (\alpha/2)$.
Take any point $E$ from the intersection of $\delta$-neighbourhoods of $\left<P_1\right>$ and $\left<P_2\right>$.
Denote by $E_1$ and $E_2$ the orthogonal projections of $E$ on $\aff P_1$ and $\aff P_2$ respectively.
Denote by $E'$ the point
of intersection $\aff P_1 \cap \aff P_2$.
We have $\angle E_1E'E + \angle E_2E'E \ge\angle E_1E'E_2$. So without loss of generality we can assume  that $\angle E_1E'E \ge \frac{1}{2} \angle E_1E'E_2$. Therefore $EE' = \frac{EE_1}{\sin \angle EE'E_1} < \frac{\delta}{\sin  (\alpha/2)} = \epsilon$. QED

\smallskip

\smallskip
{\bf Proof that $(LVKF_{k,3})\Rightarrow (LTT_{3k-1,3})$.}  Suppose to the contrary that the statement  $(LTT_{3k-1,3})$  is false for points $A_1,\ldots,A_{6k+1}$  in $\mathbb{R}^{3k-1} = \mathbb{R}^{3k-1} \times 0 \subset \mathbb{R}^{3k}$.

Without loss of generality assume that $A_1$ is a vertex of  $ \left< A_1,\ldots,A_{6k+1} \right> $.

\smallskip
{\it Construction of $6k+5$ points to apply $(LVKF_{k,3})$.}
We can assume that the points $A_1,\ldots,A_{6k+1}$ are in general position.

Consider all distances from any of these $6k+1$ points to a simplex formed by some other of these $6k+1$ points.
For any three pairwise disjoint sets $P_1,P_2,P_3$  of these $6k+1$ points consider the positive distance from  $\left<P_1\right>$ to $\left<P_2\right>\cap \left<P_3\right>$.
Choose $ \epsilon> 0 $ such that $10\epsilon$ is smaller than the smallest of all the considered numbers.

Denote by $$\pi:\mathbb{R}^{3k}\to\mathbb{R}^{3k-1} \times 0$$ the orthogonal projection.
Denote by $ \pi_K:\mathbb{R}^{3k}-\{K\}\to\mathbb{R}^{3k-1} \times 0 $ the central projection from a point $K \notin \mathbb{R}^{3k-1} \times 0$.

 Take  $\delta\in (0,\epsilon)$ given by Lemma 4 for the set $\{A_1,\ldots,A_{6k+1}\}$ and the number $\epsilon>0$.
Consider points  $M_j:=(A_1,m_j)$, $j=1,2,3,4$, where the numbers $m_j$ are defined recursively.
 Let $m_1:=1$.
If  $m_1,\ldots,m_{j-1}$ are defined, then let $m_j$ be so big that
$|\pi_{M_j}X,\pi X|<\delta$ for every $X\in M_1,\ldots, M_{j-1},A_1,\ldots,A_{6k+1}$.

 Let us prove that there exist pairwise disjoint sets  $$\Delta_1, \Delta_2, \Delta_3\subset \{M_1, M_2, M_3, M_4, A_1, A_2, ..., A_{6k+1}\}$$ such that   $\left<\Delta_1\right>\cap \left<\Delta_2\right>\cap \left<\Delta_3\right>\neq\varnothing$ and $\{M_1, M_2, M_3, M_4, A_1\}\not\subset\Delta_1\cup \Delta_2\cup \Delta_3$. Assume the contrary, then for any neighborhood of the point   $A_1$ in $\R^{3k-1}$ there exists a point  $A'_1$ from this neighborhood such that for any three pairwise disjoint subsets   $S_1, S_2, S_3 \subset \{M_1, M_2, M_3, M_4, A'_1, A_2, ..., A_{6k+1}\}$ such that  $\left<S_1\right>\cap \left<S_2\right>\cap \left<S_3\right>\neq\varnothing$ we obtain  $\{M_1, M_2, M_3, M_4, A'_1\}\subset S_1\cup S_2\cup S_3$ and $A'_1\in\left<A_1,\ldots, A_{6k+1}\right>$.
 Apply ($LVKF_{k,3}$) to the set  $\{M_1, M_2, M_3, M_4, A'_1, A_2, ..., A_{6k+1}\}$, we get pairwise disjoint sets $S_1, S_2, S_3$ such that $|S_1|=|S_2|=|S_3|=2k+1$.
Since  $\{M_1, M_2, M_3, M_4, A'_1\}\subset S_1\cup S_2\cup S_3$, one of the sets  $S_1, S_2, S_3$ contains a point  $A'_1$.
Without loss of generality, we assume that  $A'_1\in S_1$.
Since $|S_1|+|S_2|+|S_3|=6k+3$, we get some point from  $A_2, \ldots, A_{6k+1}$, which not belong in any of the sets $S_1, S_2, S_3$.
Without loss of generality, we assume that  $A_{6k+1}\notin S_1\cup S_2\cup S_3$.
Then $S_1, S_2, S_3 \subset \{M_1, M_2, M_3, M_4, A_1, A'_1, A_2, ..., A_{6k}\}$ and $A_1\notin S_1, S_2, S_3$.
Therefore, among the sets  $S_1\cap(A_1\times\R), S_2\cap(A_1\times\R), S_3\cap(A_1\times\R)$ there are at least 2 one-point ones.
Then  $\left<S_1\right>\cap \left<S_2\right>\cap \left<S_3\right>\cap(A_1\times\R)=\varnothing$.
Therefore $M_2\notin\left<S_1\right>\cap \left<S_2\right>\cap \left<S_3\right>$.
Apply Lemma 1 to $\{M_1, M_2, M_3, M_4, A_1, A'_1, A_2, ..., A_{6k}\}$, we get the sets  $S'_1, S'_2, S'_3$.
Then $|\{A_1, M_1, M_2, M_3, M_4\}\cap(S'_1\cup S'_2\cup S'_3)|\leq 3$ и $A_1\notin S'_1\cup S'_2\cup S'_3$.
This contradicts the fact that for any three pairwise disjoint subsets  $S_1, S_2, S_3 \subset \{M_1, M_2, M_3, M_4, A'_1, A_2, ..., A_{6k+1}\}$ such that $\left<S_1\right>\cap \left<S_2\right>\cap \left<S_3\right>\neq\varnothing$ we obtain  $\{M_1, M_2, M_3, M_4, A'_1\}\subset S_1\cup S_2\cup S_3$.

 We have proved that there exist pairwise disjoint sets $$\Delta_1, \Delta_2, \Delta_3\subset \{M_1, M_2, M_3, M_4, A_1, A_2, ..., A_{6k+1}\}$$ such that  $\left<\Delta_1\right>\cap \left<\Delta_2\right>\cap \left<\Delta_3\right>\neq\varnothing$ and $\{M_1, M_2, M_3, M_4, A_1\}\not\subset\Delta_1\cup \Delta_2\cup \Delta_3$.
 For this $\Delta_1, \Delta_2, \Delta_3$ we have $M_2\notin \left<\Delta_1\right>\cap \left<\Delta_2\right>\cap \left<\Delta_3\right>$.
Apply Lemma 1 to the set $\{M_1, M_2, M_3, M_4, A_1, A_2, ..., A_{6k+1}\}$.
We get the sets $\Delta'_1, \Delta'_2, \Delta'_3$ such that $\left<\Delta'_1\right>\cap \left<\Delta'_2\right>\cap \left<\Delta'_3\right> \neq \varnothing$ and there are 2 points from  $A_1, M_1, M_2, M_3, M_4$ such that do not belong to the sets $\Delta'_1, \Delta'_2, \Delta'_3$.

\smallskip

  Denote by $Z$ any point from  $\left<\Delta'_1\right>\cap \left<\Delta'_2\right>\cap \left<\Delta'_3\right>$.

{\it Case 1: $\Delta'_h\subset\{A_1,\ldots,A_{6k+1}\}$ for some $h$.}
Then $\left<\Delta'_h\right>\subset\mathbb{R}^{3k-1} \times 0$.
Hence $ Z\in\mathbb{R}^{3k-1} \times 0$.
Take $ i=1,2,3$.
Since the points $M_1, M_2, M_3, M_4$ belong to the same half-space, $\left<\Delta'_i\right> \cap (\mathbb{R}^{3k-1} \times 0) = \left< \Delta'_i \cap  (\mathbb{R}^{3k-1} \times 0) \right>$.
Then the set $Z \in \left<\Delta'_i \cap (\mathbb{R}^{3k-1} \times 0)\right>$ for $ i=1,2,3$.
Therefore sets $\Delta'_i \cap (\mathbb{R}^{3k-1} \times 0)$ for $ i=1,2,3$ constitute the required partition.

\smallskip

{\it Case 2:  none of $\Delta'_1, \Delta'_2, \Delta'_3$ is contained in $\{A_1,\ldots,A_{6k+1}\}$.}  Hence each set $\Delta'_i$  has only one vertex among the points $\{M_1, M_2, M_3, M_4\}$ and $A_1\notin\Delta'_i$ for $i=1,2,3$. Take the highest vertex of each $\Delta'_i$.
Denote these vertices by $W_1, W_2, W_3$  in the increasing order along $0\times\R$.
Without loss of generality $W_i\in \Delta'_i$  for each $i=1,2,3$.
Denote by $\Delta_i^-$ the image of  $\Delta'_i-W_i$ under the central projection $\pi_{W_i}$ from $W_i$ to $\mathbb{R}^{3k-1} \times 0$. Then $\Delta_i^-=\Delta_i\cap\R^{3k-1}\times 0$. Therefore $A_1\notin \Delta_i^-$ for $i=2,3$. Since $\pi\Delta'_1=A_1\cup(\Delta'_1\cap \mathbb{R}^{3k-1}\times 0)$, it follows that  pairwise intersections of $\pi\Delta'_1, \Delta_2^-, \Delta_3^-$ are empty.

Since $Z \in \left<\Delta'_i\right>-W_i$, it follows that $ \pi_{W_i}Z \in \left<\Delta_i^-\right>$ for $i=2,3$.
Hence
$$|\pi Z,\left<\Delta_i^-\right>|<| \pi Z, \pi_{W_i}Z |<\delta.$$
The last inequality is true by definition of $W_i$ and $M_j$.
Then by Lemma 1 the point $\pi Z$ belong to the $\epsilon$-neighbourhood of the intersection $\left<\Delta_3^-\right> \cap \left<\Delta_2^-\right>$.
Since $ Z \in \left<\Delta'_1\right>$, it follows that  $\pi Z \in \left<\pi\Delta'_1\right>$.
Hence  $|\left<\pi\Delta'_1\right>,\left<\Delta_3^-\right> \cap \left<\Delta_2^-\right>| < \epsilon$.
Then by definition of $\epsilon$ a convex hulls $\left<\pi\Delta'_1\right>,\left<\Delta_2^-\right>, \left<\Delta_3^-\right>$ have a common point.
Therefore sets $\pi\Delta'_1, \Delta_2^-, \Delta_3^-$ form a required partition.


\begin{thebibliography}{RSS95}




\bibitem[BBZ]{BBZ} \emph{I. Barany, P. V. M. Blagojevic and G. M. Ziegler.}, Tverberg's Theorem at 50: Extensions and Counterexamples, Notices of the AMS, 63:7 (2016), 732-739. http://www.ams.org/journals/notices/201607

\bibitem[BFZ14]{BFZ14} \emph{P. V. M. Blagojevic, F. Frick, and G. M. Ziegler}, Tverberg plus constraints, Bull. Lond. Math. Soc. 46:5 (2014), 953-967, arXiv: 1401.0690.

\bibitem[BZ]{BZ} \emph{P. V. M. Blagojevic and G. M. Ziegler.}, Beyond the Borsuk-Ulam theorem: The topological Tverberg story, arXiv:1605.07321v2

\bibitem[Fr15]{Fr15} \emph{F. Frick}, Counterexamples to the topological Tverberg conjecture,
arXiv:1502.00947.

\bibitem[Gr10]{Gr10} \emph{M. Gromov}, Singularities, expanders, and topology of maps. Part 2: from combinatoricsto topology via algebraic isoperimetry., Geometric and Functional Analysis, 20:2 (2010), 445–446.

\bibitem[Sk17]{Sk17} \emph{A. Skopenkov}, On van Kampen-Flores, Conway-Gordon-Sachs and Radon theorems, published as \S4
Russian Math. Surveys, 73:2 (2018), 323-353
arxiv:1704.00300

\bibitem[Sk18]{Sk18} \emph{A. Skopenkov}, A user's guide to the topological Tverberg
conjecture,  Russian Math. Surveys,  73:2(2018) 323-353,
arxiv:1605.05141.

\bibitem[Tv66]{Tv66} \emph{H. Tverberg} (1966), A Generalization of Radon's Theorem. Journal of the London Mathematical Society, s1-41: 123-128. doi:10.1112/jlms/s1-41.1.123


\end{thebibliography}
\end{document}